\documentclass[12pt]{article}
\usepackage{latexsym,amsfonts,epsf}

\begin{document}

% theorem environments
\newtheorem{thm}{Theorem}[section] 
\newtheorem{lem}[thm]{Lemma}
\newtheorem{cor}[thm]{Corollary} 
\newtheorem{deff}[thm]{Definition}
\newtheorem{examm}[thm]{Example}
\newtheorem{prooff}{Proof:}
\newtheorem{axiom}{Axiom}

\newenvironment{define}{\begin{deff}\normalfont}{\end{deff}}	%unitalic
\newenvironment{example}{\begin{examm}\normalfont}{\end{examm}}	%unitalic
\newenvironment{proof}{\begin{prooff}\normalfont}{\end{prooff}}	%unitalic
\renewcommand{\theprooff}{}					%no numbers

% other definitions
\newcommand{\integers}{\ensuremath{\mathbb{Z}}}
\newcommand{\reals}{\ensuremath{\mathbb{R}}}
\newcommand{\relators}{\ensuremath{\mathcal{R}}}
\newcommand{\burnside}{\ensuremath{\mathcal{B}}}
\newcommand{\cayley}{\ensuremath{\mathcal{C}}}
\newcommand{\poincare}{\ensuremath{\mathcal{P}}}
\newcommand{\windnum}{\mbox{wind.num.}} 	% for lack of a better idea
\newcommand{\str}{\mbox{str}}
\newcommand{\cir}{\mbox{cir}}
\newcommand{\rank}{\mbox{rank}}

\title{The Burnside Groups\\ and Small Cancellation Theory}
\author{Jonathan P. McCammond} 
%\address{Department of Mathematics\\ 
%	Texas A \& M University \\ 
%	College Station, TX 77843}
%\email{jon.mccammond@math.tamu.edu}
%\subjclass{20F50}
%\keywords{Burnside groups, small cancellation theory} 
\date{\today}

\maketitle

\begin{abstract}

In a pair of recent articles\footnote{As of 1 Dec 1997, one article
has been accepted for publication and the second will soon be
available as a preprint.  Until the second article has been
distributed, all of the results which refer to the Burnside groups in
particular should be considered preliminary announcements only.  The
exponent of 1260 is especially provisional.}, the author develops a
general version of small cancellation theory applicable in higher
dimensions (\cite{mccammond-gsct}), and then applies this theory to
the Burnside groups of sufficiently large exponent
(\cite{mccammond-burnside}).  More specifically, these articles prove
that the free Burnside groups of exponent $n \geq 1260$ are infinite
groups which have a decidable word problem.  The structure of the
finite subgroups of the free Burnside groups is calculated from the
automorphism groups of the general relators used in the
presentation. The present article gives a brief introduction to the
methods and techniques involved in the proof.  Many of the ideas
originate with the recent work of A.~Yu.~Ol'shanskii and S.~V.~Ivanov.
Some familiarity with their work on the Burnside problem will be
assumed.

\end{abstract}

\section{Introduction}\label{sect-introduction}

Ninety-five years ago William Burnside asked whether every finitely
generated group of finite exponent must be finite \cite{burnside}. 
Since then, significant progress has been made in the study of the free
Burnside groups, much of it coming in recent years. In this first
section the recent history of their study will be briefly
reviewed so that the results announced below can be viewed in their proper
context.  The second section introduces many of the concepts from which
the general theory is constructed.  The third and the fourth sections
will describe the results of general small cancellation theory and the
results on the Burnside groups in greater detail.

\subsection{Historical Context}

% Adian and Novikov

The existence of exponents for which Burnside's conjecture is false was
first demonstrated by Novikov and Adian in a series of articles
published in 1968~\cite{adiannovikov}. More specifically, for $n \geq
4381$ and odd, they proved that the free Burnside group of exponent $n$
is infinite, that these groups have decidable word problems, and that
their only finite subgroups are cyclic.  In a book-length reworking of
the results \cite{adian-english}, Adian was able to extend coverage of
the proof to include all cases with $n \geq 665$ and odd.  Both
versions of the proof are several hundred pages in length.

% Ol'shanskii

In 1982, A.~Yu.~Ol'shanskii published an alternative proof that the
Burnside groups of exponent $n > 10^{10}$ and odd are infinite, that
they have a decidable word problem, and that their only finite
subgroups are cyclic \cite{olshanskii}.  Ol'shanskii used small
cancellation theory and especially van Kampen diagrams to create a
significantly shorter proof of Adian and Novikov's results.  By
introducing more geometry, the proof was reduced to a mere thirty-page
article.

% Ivanov

In 1994, S. V. Ivanov published a substantial article on the Burnside
groups in which he showed that when $n$ is at least $2^{48}$ and either
odd or divisible by $2^9$, then the free Burnside group of exponent $n$
is infinite~\cite{ivanov}.  In addition, he showed that these groups
have a decidable word problem and he effectively described their finite
subgroups.  Since every Burnside group with exponent at least $2^{57}$
has one of these groups as a homomorphic image, Ivanov's proof shows
that Burnside's original conjecture can be answered in the negative for
almost all exponents $n$.  

% Lysionok

Finally, in 1996 I.~Lysionok published a proof that the Burnside
groups of exponent $n = 16k \geq 8000$ are infinite \cite{lysionok}.
Combined with the earlier proofs regarding Burnside groups of large
odd exponent, this work also shows that almost all of the Burnside
groups are infinite.

% summary

As can be seen from the above review, the study of an arbitrary Burnside
group of large exponent has proceeded indirectly through the study of
the homomorphic images onto Burnside groups whose exponent is either
large and odd or divisible by a large power of 2.  In particular, the
restriction that $n$ be either large and odd or divisible by a large
power of 2 has always been needed to show the decidability of the word
problem, the structure of the finite subgroups, and other detailed
structural information.  The results described below remove these
restrictions.

\subsection{Announced Results}

% McCammond

The author is currently in the process of publishing and preparing to
publish a pair of articles which first develop a general version of small
cancellation theory applicable in higher dimensions
(\cite{mccammond-gsct}), and then apply this theory to the Burnside
groups of sufficiently large exponent (\cite{mccammond-burnside}). In
combination these articles prove that the free Burnside groups of
exponent $n \geq 1260$ are infinite groups which have a decidable word
problem.  The structure of the finite subgroups of these groups is
calculated from the automorphism groups of the general relators used in
the presentation. The present article gives a brief introduction to the
methods and techniques involved in the proof. The definition of a
general relator, the axioms of general small cancellation theory, and a
description of the inductive construction involved for the Burnside
groups are included.  Many of the ideas involved originate with the
recent work of A.~Yu.~Ol'shanskii and S.~V.~Ivanov.  Some familiarity
with their work on the Burnside problem will be assumed. 

The theory developed by the author in \cite{mccammond-gsct} is a
generalized version of small cancellation theory which is applicable to
specific types of high-dimensional simplicial complexes. The usual
results on small cancellation groups are shown to hold in this new
setting with only slight modifications.  The major results of this
theory are summarized below in Theorem~\ref{mainthm1}.  The
definitions of the unfamiliar terms will be given later in the article.

\renewcommand{\thethm}{\Alph{thm}}   % to label these Theorem A and B

\begin{thm}\label{mainthm1}

If $G = \langle A | \relators \rangle$ is a general small cancellation
presentation with $\alpha \leq \frac{1}{12}$, then the word and
conjugacy problems for $G$ are decidable, the Cayley graph is
constructible, the  Cayley category of the presentation is
contractible, and $G$ is the direct limit of hyperbolic groups.  If the
presentation satisfies a few additional hypotheses, then every finite
subgroup of $G$ is a subgroup of the automorphism group of some general
relator in $\relators$.  

\end{thm}

The generalized small cancellation theory developed in
\cite{mccammond-gsct} overlaps significantly with the general theories
developed by Rips (\cite{rips}) and Ol'shanskii (\cite{olshanskii}).
In fact, the theory described below attempts to provide an underlying
geometric object whose existence explains the success of these earlier
constructions.

The main results on the Burnside groups which are proved in
\cite{mccammond-burnside} are summarized below in
Theorem~\ref{mainthm2}. The notation $\burnside(m,n)$ refers to the
$m$-generated free Burnside group of exponent $n$.

\begin{thm}\label{mainthm2}

For any $m>1$, and any $n \geq 1260$, the Burnside group
$\burnside(m,n)$ possesses a general small cancellation presentation
which satisfies all of the hypotheses and conclusions of
Theorem~\ref{mainthm1}.  In addition, every finite subgroup of the
group $\burnside(m,n)$ is contained in a direct product of a dihedral
of order $2n$ with a finite number of dihedral $2$-groups whose exponent
divides $n$.  As a consequence, all of the groups $\burnside(m,n)$ are
infinite.

\end{thm}

\renewcommand{\thethm}{\thesection.\arabic{thm}}    % return to normal

The approach used in these articles has succeeded in significantly
lowering the previous bounds beyond which all of the Burnside groups
are known to be infinite, and also in proving for the first time that
the word problem for all Burnside groups of sufficiently large exponent
is decidable.   Finally, in contrast with earlier work on the Burnside
groups, no distinction is made between Burnside groups with even or odd
exponents.

\section{Preliminaries}\label{sect-prelim}

The generalization of small cancellation theory alluded to above  will
depend on an expansion of the usual definition of a relator.  To help
motivate the revised definition, we will begin with a few familiar
examples, followed by the technical definitions.

\subsection{First Examples}

The following examples are groups with familiar presentations and
well-known properties which can be subsumed under the rubric of general
small cancellation theory.  The explanation of the exact definitions will
be postponed until after the examples have been described.

\begin{figure}
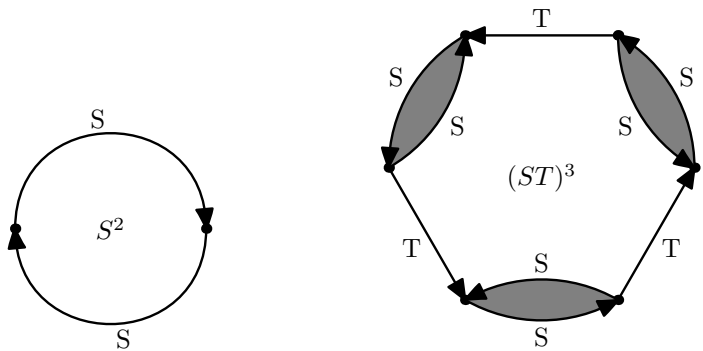

	\begin{center}~
		\epsffile{modular.1}
		\hspace{5em}
		\epsffile{modular.2}
	\end{center}
	\caption{General relators for the modular group \label{fig-modular}}
\end{figure}

\begin{example}\label{exam-modular}

{\bf [The Modular Group]} Let 

\[ 	
	S = \left(
	\begin{array}{cc}
		0 & -1 \\ 1 & 0
	\end{array}
	\right) \textrm{ and }
	T = \left(
	\begin{array}{cc}
		1 & 1 \\ 0 & 1
	\end{array}
	\right)
\]

\noindent Then $S^2 = -I, (ST)^3 = -I$.  The modular group
$\textrm{PSL}_2(\integers) = \textrm{SL}_2(\integers) / \pm 1$ is the
free product of the cyclic groups of order 2 and 3 which are generated
by $S$ and $ST$.  In particular, $\textrm{PSL}_2(\integers)$ is
generated by $S$ and $T$, and it has the following presentation:

\[
	\textrm{PSL}_2(\integers) = \langle S,T | S^2 = (ST)^3 = 1 \rangle
\]

\noindent This particular presentation does not satisfy any of the usual
small cancellation conditions such as $C(6)$ or $C'(\frac{1}{6})$ since
the letter $S$ is a piece and it represents one-half of the boundary of
the first relator.  The general small cancellation theory under
discussion would replace this presentation with a
`general presentation' which still has $S$ and $T$ as generators, but the
relators will be altered to appear as in Figure~\ref{fig-modular}.

\end{example}

\begin{figure}
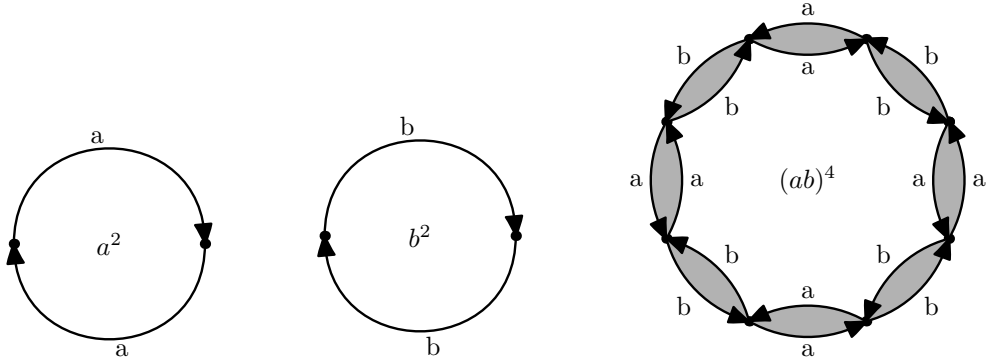

	\begin{center}~
		\epsffile{dihedral.1}
		\hspace{3em}
		\epsffile{dihedral.2}
		\hspace{3em}
		\epsffile{dihedral.3}
	\end{center}
	\caption{General relators for the dihedral group of order 8
	\label{fig-dihedrals}} 
\end{figure}

\begin{example}\label{exam-dihedrals}

{\bf [Dihedral Groups]}  Consider the standard presentation of the
dihedral group of order $2n$:

\[ 
	\textrm{D}_{2n} = \langle a,b | a^2 = b^2 = (ab)^n = 1 \rangle
\]

\noindent This presentation also has a modification which satisfies
some general small cancellation conditions.  The resulting general
relators are given in Figure~\ref{fig-dihedrals} Moreover, it is easy
to show that no presentation of a dihedral is ever a small
cancellation presentation in the traditional sense, since one of the
results of the traditional theory is that all of the finite subgroups
of a small cancellation group are cyclic.

\end{example}

\begin{example}\label{exam-coxeter}

{\bf [Coxeter Groups]}  Coxeter groups can be seen as a generalization
of Example~\ref{exam-dihedrals}.  Let $M$ be a symmetric $n \times n$
matrix whose entries lie in the set $\integers^+ \cup \infty$.  Assume
in addition that $m_{ii} = 1$ and $m_{ij} > 1$ for all $i \neq j$.  The
Coxeter group based on $M$ is given by the presentation

\[
	G = \langle a_1, a_2, \ldots , a_n | (a_ia_j)^{m_{ij}} = 1 \rangle
\]

\noindent where no relation is added in the case $m_{ij} = \infty$.  If
$m_{ij}$ is never equal to 2, then the Coxeter group is said to be of
large type.  If $m_{ij}$ is never equal to 2 or 3, then the Coxeter group
is said to be of extra-large type.  Once the relators $(a_ia_j)^{m_{ij}}$
are modified as in Example~\ref{exam-dihedrals}, the resulting Coxeter
group satisfies a version of $C(4)$, a Coxeter group of large type satisfies
a version of $C(6)$, and a Coxeter group of extra-large type satisfies
a version of $C(8)$.

\end{example}

\subsection{General Relators}

The most important concept involved in general small cancellation
theory is that of a general relator.  Every general relator will be
assigned a height.  Traditional relators are examples of height 2.
The examples given in the previous section are general relators of
height 3.  Before the precise technical definition is given, the
examples from the previous section will be discussed in greater
detail, and we will present a single example of a general relator of
height 4.

\begin{figure}
	\begin{center}~
		\epsffile{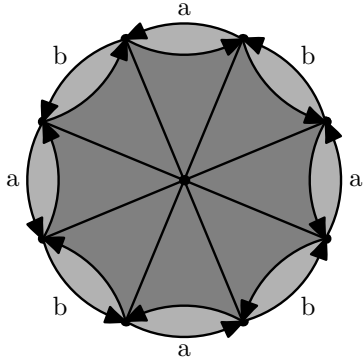}
	\end{center}
	\caption{A more precise picture of a general relator\label{fig-solid}} 
\end{figure}

\begin{define}\label{def-boundary}

{\bf [Boundaries of General Relators]} In the same way that a relator
in the traditional sense is often viewed as a labeled disk rather than
a labeled circle, the structures shown in Figure~\ref{fig-modular} and
Figure~\ref{fig-dihedrals} are, technically speaking, merely the
boundaries of general relators and not the general relators
themselves.  A general relator will be the topological cone over its
roughly circular boundary.  Thus, the general relator corresponding to
$(ab)^4$ (whose boundary is shown in Figure~\ref{fig-dihedrals}) will
be composed of eight solid circular cones which are adjoined along
their lateral sides.  This particular general relator is shown
schematically in Figure~\ref{fig-solid}.  Similarly the general
relator corresponding to $(ST)^3$ (whose boundary is shown in
Figure~\ref{fig-modular}) will consist of three triangles and three
solid circular cones.

\end{define}

\begin{define}\label{def-height}

{\bf [Heights]} In a traditional relator, the vertices, the open edges
and the interior of the disk can be thought of as being represented by
their barycenters.  Moreover, these barycenters can be assigned
heights based on whether they come from a vertex (height 0), an edge
(height 1), or a disk (height 2).  In a similar way, the apex of the
cone over the dihedral structure in Figure~\ref{fig-solid} can be
assigned a height of 3 since it contains a disk of height 2 in its
boundary.  Once these barycenters have been assigned heights, there is
a notion of a 1-, 2-, and 3-skeleton which is the union of all of the
pieces assigned a height of at most 1, at most 2, or at most 3.  The
rightmost picture in Figure~\ref{fig-dihedrals} is thus the 2-skeleton
of the general relator, as well as its boundary.

\end{define}

\begin{figure}
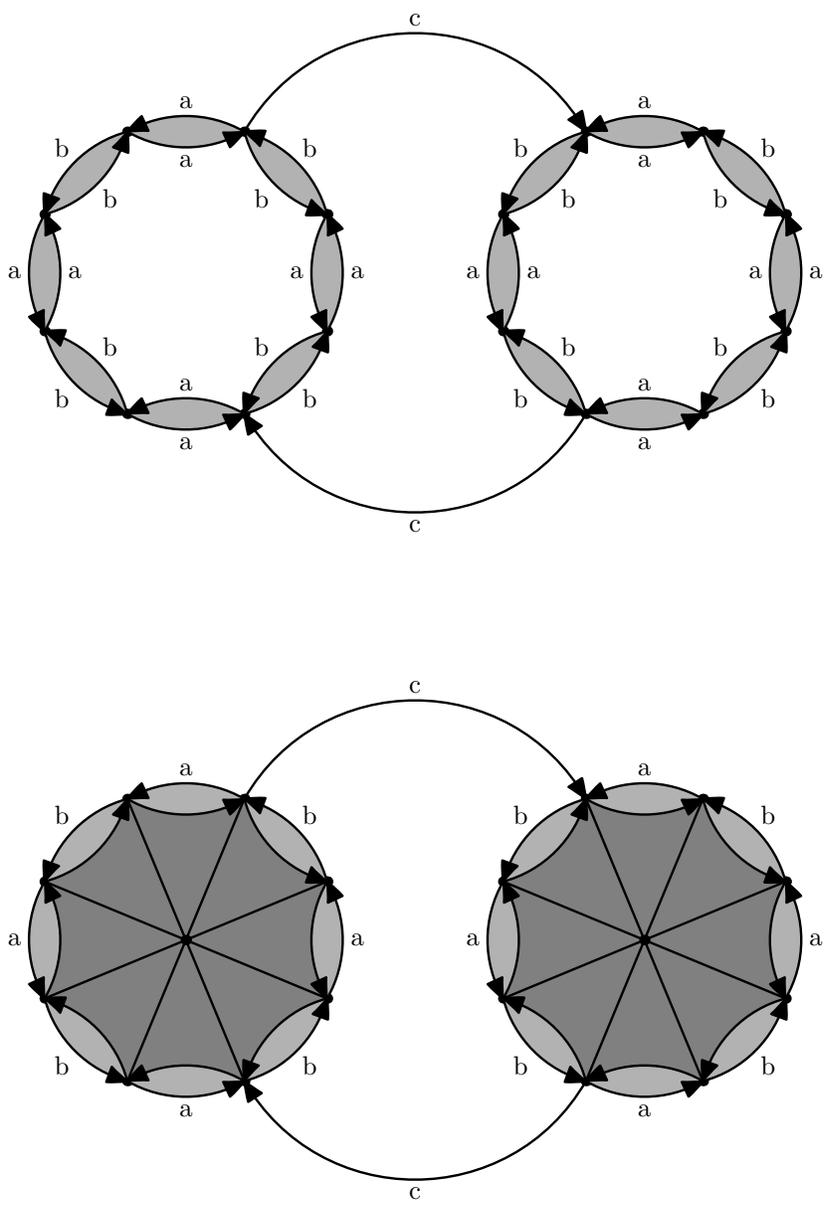

	\begin{center}~
		\epsffile{highrank.1}
	\end{center}
	\vspace{.4in}
	\begin{center}~
		\epsffile{highrank.2}
	\end{center}
	\caption{The 2- and 3-skeletons of a general relator of height 4
	\label{fig-highrank}} 
\end{figure}

\begin{example}\label{exam-higherheight}

{\bf [A General Relator of Height 4]} The following explicit example
will be used to illustrate several key concepts which are precisely
defined below.  Let $G = \langle a,b,c | a^2, b^2, (ab)^4, (abcb)^2
\rangle$.  A modification of this presentation which uses general
relators begins with the three (general) relators whose boundaries are
shown in Figure~\ref{fig-dihedrals}.  The fourth relator, $(abcb)^2$,
is modified as shown in Figure~\ref{fig-highrank}.  It is obtained
from the cycle labeled $(abcb)^2$ by attaching the general relator for
$(ab)^4$ to each of the paths labeled $bab$.  The top figure shows the
2-skeleton of this relator, while the bottom figure shows its
3-skeleton.  As described above, the general relator itself will be
the topological cone over this roughly circular boundary, and the apex
of the cone will be assigned a height of 4 since it contains a point
of height 3 in its boundary.

\end{example}

Examples of general relators of higher height can be obtained as
follows:  Given a general relator $R$ of height $k$, introduce a new
generator $d$ and attach an arrow labeled $d$ to any two distinct
vertices in $R$.  The resulting structure qualifies as the boundary of
a general relator, and since it contains $R$ in its boundary it has
height $k+1$.  Variations on this theme can create presentations which
contain general relators of arbitrary height.

As can be seen from these examples, a general relator is no longer
simply a word or a cycle, but rather a particular type of simplicial
complex with a labeling on its 1-skeleton.  The detailed structure of
a general relator is best captured by the partially ordered set of its
barycenters.  In each of the examples given so far we can associate a
partially ordered set to each general relator.  The elements are the
barycenters, and the ordering is defined by setting $p < q$ if and
only if $p$ is contained in the boundary of $q$.  Perhaps more
surprising, this partially ordered set, or poset, contains all of the
information needed to reconstruct the complex.  (The procedure for
creating a simplicial complex from a poset is due to D. Quillen
\cite{quillen}).  For an arbitrary partially ordered set the
corresponding simplicial complex is derived by taking the finite
chains as the simplices.  For a traditional relator, the poset has
only 3 levels: the vertices are the elements of height 0, the open
edges are the elements of height 1, and the single open 2-cell is the
unique element of height 2.  The simplicial complex corresponding to
this poset is a simplicial subdivision of the original 2-cell.  The
$k$-skeleton of such a poset will refer to the geometric realization
of the sub-poset consisting only of those elements of height $k$ or
less.  Thus in the traditional case, the 1-skeleton is the cyclic
boundary and the 0-skeleton consists of the vertices alone.  For the
general relator shown in Figure~\ref{fig-highrank}, the poset has one
element of height 4, two elements of height 3, sixteen elements of
height 2, and sixteen elements of height 1.

We will now precisely define a general relator using these partially
ordered sets.

\begin{define}\label{def-relators}

{\bf [General Relators]} A general relator $R$ is (the simplicial
complex corresponding to) a finite poset with a unique maximum element
subject to certain additional restrictions. First, the 1-skeleton of
the poset must be a simplicial subdivision of a graph, and the
unsubdivided graph needs to be deterministically labeled by a set of
generators.  Next, for all elements $p$ of height at least 2, the
geometric realization of the ideal of all elements strictly below $p$
must be homotopically equivalent to a circle.  In other words, the
simplicial complex corresponding to this ideal should be homotopically
equivalent to the space $S^1$.  If the unique maximum element has
height $k$ then the general relator $R$ is said to be of height $k$,
and the $(k-1)$-skeleton of $R$ is called its boundary. General
relators were introduced by the author in \cite{mccammond-gsct}.

\end{define}

Notice that the boundary of a general relator is allowed to be fairly
complicated as long as it is homotopically equivalent to the unit
circle.  Although a general relator can contain a 1-skeleton which can
be significantly more complicated than is allowed in the traditional
theory, the local structure is often less important than its global
topology.  One example of a global topological concept which plays a
key role in the theory is that of the winding number of a loop.  A
winding number is definable in this context because of the homotopic
equivalence to the unit circle.  A loop in the boundary of a general
relator with winding number 1 will be called a representative of the
general relator.  Using representatives, it is possible to define more
or less traditional van Kampen diagrams over collections of general
relators by requiring that the label of every 2-cell in the planar van
Kampen diagram be the label of a representative loop in the boundary
of some general relator.

\begin{define}\label{def-representatives}

{\bf [Representatives]} A representative of a general relator $R$ is any
path in $R$ which forms a loop with winding number $\pm 1$.  There are
no other restrictions on representatives:  they can complete five
clockwise rotations followed by four counterclockwise ones, they can
self-intersect, they can even be unreduced in the free group.  The main
use to which these representatives will be put is to define a van Kampen
diagram over a set of general relators.

\end{define}

\begin{define}\label{def-vKdiagrams}

{\bf [Van Kampen Diagrams]}  A van Kampen diagram over a set $\relators$
of general relators is exactly like a van Kampen diagram over a set of
traditional relators except that the boundary cycle of each 2-cell in the
diagram is a path which maps onto a representative path in one of
the general relators in $\relators$, instead of being a relator itself.

\end{define}

In traditional small cancellation theory van Kampen diagrams are reduced
and cancellable pairs are removed until in the reduced diagram a
curvature condition forces a large portion of a relator to exist on
the boundary.  In the general theory described here a similar strategy is
pursued, but one situation arises which does not need to
be considered in the traditional case.

\begin{example}\label{exam-self}

{\bf [Self-bordering Cells]}  The annular diagram shown in
Figure~\ref{fig-self} is a van Kampen diagram over the presentation
given in Figure~\ref{fig-dihedrals}.  More explicitly, starting at the
rightmost vertex, and reading counterclockwise around the outer circle,
the boundary of the 2-cell in the diagram reads the word $W =
ababa^{-1}b^{-1}a^{-1}b^{-1}$, and this is a path which is a
representative of the rightmost relator in Figure~\ref{fig-dihedrals}.
The difficulty that this diagram causes is that the word $bab$ is long
relative to the length of the boundary.  In the traditional theory, an
overlap longer than a piece is prohibited by the fact that no
(nontrivial) word in the free group is conjugate to its own inverse. 
The word $W$, however, is conjugate to its own inverse, as
Figure~\ref{fig-dihedrals} demonstrates.

\begin{figure}
	\begin{center}~
		\epsffile{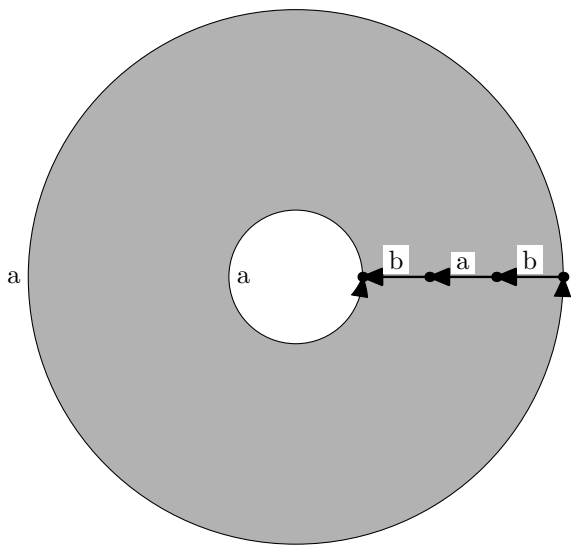}
	\end{center}
	\caption{A 2-cell with a large overlap \label{fig-self}}
\end{figure}

\end{example}

When diagrams such as this are encountered as subdiagrams during the
reduction of a van Kampen diagram, there must be some way to remove
them since the goal of the reductions is to produce a reduced diagram
with no long internal arcs.  For Example~\ref{exam-self} the reduction
is easy since the entire annular subdiagram can be removed and the
inner and the outer loops labeled $a$ can be identified.  (The general
situation is more complicated but those details will be glossed over
here.)

There is another way to explain the existence of self-bordering cells
in the general case.  Notice that corresponding to a self-bordering
2-cell which represents a general relator $R$ there is a
label-preserving automorphism of the boundary of $R$ which reverses
the orientation of its boundary.  In particular, the self-bordering
2-cell in Figure~\ref{fig-self} corresponds to the automorphism of
Figure~\ref{fig-solid} which rotates the left- and right-hand $a^2$
loops by 180 degrees.  This automorphism will preserve the labels on
the edges but it will reverse the orientation of the circular
boundary.  Traditional relators cannot have automorphisms of this
type.

Although general relators can have automorphisms which are more
complicated than traditional relators, the structure of the possible
automorphism groups is fairly restricted by group theoretic standards.
Before the relevant theorem is quoted, however, it is necessary to make
plain the meaning of two of the expressions used in its statement.
First, the phrase `closed under subcones' is just another way of
stipulating that a general relator which is contained in the boundary
of a general relator in the set $\relators$ must also be contained in
$\relators$.  Thus a set of general relators which contains the one
whose boundary is shown in Figure~\ref{fig-highrank} must also contain
the general relator shown in Figure~\ref{fig-solid}.

Next, a crucial cone in the boundary of a general relator $R$ is,
loosely speaking, an edge or a general relator whose removal from the
boundary of $R$ makes it impossible to complete a representative.  In
a traditional relator all of its edges are crucial.  In the general
case, there are some additional restrictions on the way in which the
crucial cone is situated in the boundary. See \cite{mccammond-gsct}
for details.

\begin{thm}[\cite{mccammond-gsct}]\label{thm-autR}  % Lemma 5.26

Let $\relators$ be a set of general relators closed under subcones and
suppose that all general relators in $\relators$ have at least one
crucial cone in their boundary.  If $H$ is a subgroup of the
automorphism group of a general relator $R \in \relators$, then for
some $r$ there is a group homomorphism from $H$ to $\textrm{D}_{2r}$
whose kernel is a 2-group.  In particular, $H$ is isomorphic to a
cyclic or a dihedral group extended by a 2-group.

\end{thm}

The main advantage which general relators present over the traditional
type is that certain symmetries which exist in the group can be
captured by the relations.  This obviates the need to break these
symmetries by making arbitrary choices, such as selecting one
representative loop of length 8 in the final relator in
Figure~\ref{fig-dihedrals} and not any of the others.  When these
symmetries are incorporated into the general relators themselves, the
structure of the group beyond the relator level becomes more visible.

\subsection{Length and Width}

Because general relators are allowed to have boundaries which are
`thick', the portion of the boundary traversed by a path requires a bit
of exposition.  Along the way, a length and a width will be prescribed
for each general relator.  A single example will be sufficient to
illustrate the nuances of these definitions.

\begin{example}\label{exam-moebius}

{\bf [Twisted Boundaries]}  Consider the presentation

\[
	G = \langle a,b,c | (ab)^3 = (bc)^3 = (ca)^3 = W = 1 \rangle
\]

\noindent  where $W = babcba^{-1}c^{-1}$.  If the relator $W$ had been
excluded, the presentation would satisfy $C(6)$.  The relator $W$ can
be altered to form a general relator by attaching a relator to $W$
whenever the cycle $W$ contains a sizable portion of this other
relator. This process quickly stops at a structure which is a
M\"{o}bius strip. The result is shown in Figure~\ref{fig-moebius}
except that the arrow on the far left labeled $a$ must be glued (in
an orientation-preserving way) to the arrow on the far right labeled
$a$.   The word $W$ can be traced in the diagram from the upper left 
to the lower righthand corner.  Notice also that this path becomes a
loop under the gluing operation.

The resulting structure satisfies all of the conditions necessary to
qualify as the boundary of a general relator.  Let $R$ be the name
of the general relator which results.  The general relator $R$ is
technically the topological cone over the M\"{o}bius strip, the
M\"{o}bius strip being merely the boundary of $R$.  The twisting of
the boundary which occurs in this example occurs in a great many
general small cancellation presentations, including those constructed
for the Burnside groups.

\begin{figure}
	\begin{center}~
		\epsffile{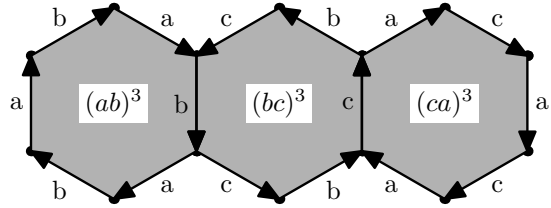}
	\end{center}
	\caption{A construction based on the word $W$ 
	\label{fig-moebius}} 
\end{figure}

\end{example}

The general relator $R$ decribed above will now be used to illustrate
several key definitions. The first such concept is that of length of a
general relator.  For a traditional relator, the length of its boundary
is simply the length of the unique reduced loop of winding number 1. In
the general case the situation is more subtle.  In
Example~\ref{exam-moebius} the loop $W$ is the shortest representative
of general relator $R$ and it has length 7.  But among those loops with
winding number 2, some have length 12.  For loops with
extremely large winding number the geodesic length of the path is more
closely approximated by six times the winding number rather than by
seven.  In order to accurately estimate large geodesic distances, the
length of a general relator must be based on the smallest average
length per winding number.  Thus the length of the general relator $R$
will be defined to be 6 even though there is no representative loop
with that length.

\begin{define}\label{def-length}

{\bf [Length of a General Relator]} The length of a general relator $R$
is defined to be the smallest value of the length of a loop in $R$ divided
by its nonzero winding number.  It is easy to show that such a minimum value
must exist, and to create examples where the length of the relator is a
fraction instead of an integer.  The length of $R$ will be denoted $|R|$.
Finally, notice that for traditional relators, this procedure gives the
traditional value for the length.

\end{define}

Now that the length of the relator itself has been clarified, it is
possible to define the relative length of a path in the boundary of a
general relator.

\begin{define}\label{def-graphmetric}

{\bf[Graph Metric]}  The (normalized) graph metric for a general relator
$R$ is a function which assigns a nonnegative length to every path in
the boundary of $R$.  In particular, if $U$ is a path in $R$ whose lift
to the universal cover of the boundary of $R$ is a geodesic path, then
the length of $U$ relative to $R$, denoted $|U|_R$, is the length of
$U$ divided by the length of $R$.  If the path $U$ is not a geodesic
when it is lifted to the universal cover of the boundary of $R$, then
let $V$ be a geodesic path in the universal cover between the same
vertices as $U$ and set $|U|_R$ equal to $|V|/|R|$. Equivalently, the
normalized length of $U$ is the length of the shortest path which is
homotopic to $U$ in the boundary of $R$ relative to its endpoints,
divided by the length of $R$ defined above.

\end{define}

Consider the path $W$ in the relator $R$ from
Example~\ref{exam-moebius}. Since the length of the general relator
$R$ is 6, since the length of the word $W$ is 7, and since $W$ is a
geodesic when lifted to the universal cover of the boundary of $R$,
$|W|_R = \frac{7}{6}$.  A metric similar (and possibly identical) to
the one just described will be used to test whether a set of general
relators satisfies a general small cancelation condition.

An additional concept to be illustrated by Example~\ref{exam-moebius} is
that of width.  Let $R^\infty$ represent the universal cover of the
boundary of a general relator $R$.  Since the boundary of $R$ is by
definition topologically equivalent to a circle, it seems clear that
the universal cover should look topologically like the real line.  Thus
$R^\infty$ should have two infinite `ends', one going in the positive
direction and the other going in the negative.   These intuitions are
true and they can be made precise but we will not do so here.  See
\cite{mccammond-gsct} for further details.

The width of a general relator is the smallest integer needed to
guarantee that the infinite ends of $R^\infty$ will be disconnected
from each other upon the removal of an open ball, regardless of the
vertex used for the center.

\begin{define}\label{def-width}

{\bf [Width of a General Relator]}  The width of a general relator $R$,
denoted $\omega_R$, is the unique smallest integer needed to guarantee
that the removal of an open ball of radius $\omega_R$ centered at an
arbitrary vertex in $R^\infty$ will always disconnect the two infinite
ends. Technically, only the 1-skeleton must be disconnected; the 
higher-dimensional skeleta are ignored. Also, for convenience, an open
ball will be defined so that it always contains the vertex used as its
center.  As a result, the width of a traditional relator is 0.

\end{define}

In Example~\ref{exam-moebius}, the width of the general relator $R$ is 3.
To see this, consider the vertex directly below the label $(bc)^3$ in
Figure~\ref{fig-moebius}, and call this vertex $u$.  Assume in addition
that the figure is a portion of the universal cover, so that the row of
adjacent hexagons extends infinitely in both directions.  An open ball of
radius 1 centered at $u$ will remove only the vertex $u$ and the two open
edges incident with it.  An open ball of radius 2 will remove a total
of six edges and three vertices, but a path running along the top edge of
the structure still succeeds in connecting the ends.  With a radius of
3, the ends are disconnected, and the three connected components
which result are the positive end, the negative end, and an isolated
vertex which is directly above the vertex $u$.

\subsection{General Presentations}

A general presentation is a set of generators $A$ together with a
collection of general relators $\relators$ which are labeled by $A$.
Such a presentation will be denoted $G = \langle A | \relators
\rangle$. The group assigned to a particular general presentation is
based on the fundamental group of a topological space constructed from
the generators and relators.  For traditional presentations 
the standard construction of a 2-complex has the described group
as its fundamental group.  This complex is sometimes referred to as
Poincar\'{e}'s construction.  For a general presentation there is a
variation on this procedure. The resulting complex will be called the
Poincar\'{e} construction of the presentation.

\begin{define}\label{def-poincare}

{\bf[Poincar\'{e} Constructions]} The Poincar\'{e} construction of a
general presentation is formed by taking a single vertex, attaching an
edge for each of the generators, and then attaching the general
relators of successively higher height, in order.  The attaching of
the general relators of height 2 yields a structure that is precisely
the standard 2-complex corresponding to these generators and
relations. If the set of general relators, $\relators$, is closed
under subcones and the relators are attached in order of increasing
height, then when a relator $R$ is attached, all the edges and general
relators used to construct the boundary of $R$ have already been added
to the Poincar\'{e} construction.  In addition there exists a unique
map from the boundary of $R$ into the structure as constructed so
far. The general relator $R$ is then attached using the map from the
boundary to the developing structure as the attaching map.  After all
the relators have been added, the result is the completed Poincar\'{e}
construction.

\end{define}

\begin{example}

The idea behind this construction can be illustrated using
Example~\ref{exam-dihedrals}.  Of the three general relators in this
presentation, two are traditional relators of height 2 and the final
relator is of height 3.  The Poincar\'{e} construction for this
presentation begins by attaching the disks labeled $a^2$ and $b^2$ to
a bouquet of two circles labeled $a$ and $b$.  The construction at
this point is the standard 2-complex whose fundamental group is the
infinite dihedral group.  The final general relator is now attached
along its boundary to the construction.  Its boundary is illustrated
in Figure~\ref{fig-dihedrals}.  The relator itself is shown in
Figure~\ref{fig-solid}.  Once attached the new complex will have the
dihedral group of order 8 as its fundamental group.

\end{example}

\begin{define}\label{def-cayley}

{\bf[Cayley Categories]} The Cayley category of a presentation is
based on the universal cover of the Poincar\'{e} construction, but the
duplication which arises in the universal cover is removed.  The
details are as follows.  The universal cover of the Poincar\'{e}
construction contains a 1-skeleton which is the Cayley graph of this
group.  As in the traditional theory, the universal cover of the
Poincar\'{e} construction may contain multiple copies of a general
relator attached to this Cayley graph by functions which agree on the
1-skeletons in their boundaries.  This occurs precisely when the
general relator under consideration possesses nontrivial
automorphisms. The number of such multiplicities is governed by the
size of the automorphism group of the particular general relator
involved.  Once these multiple copies are eliminated through a
suitable identification process, the resulting structure is called the
Cayley category of the presentation.  The map from the universal cover
to the Poincar\'{e} construction factors through the Cayley category.

\end{define}

\begin{example}

If we again use the dihedral group of order 8 as an example, the
universal cover of the Poincar\'{e} construction has the Cayley
graph as its 1-skeleton, and has two disks labeled $a^2$ attached
to each $a^2$ loop, and two disks labeled $b^2$ attached to each $b^2$
loop.  If we forget the way they are mapped into the Poincar\'{e}
construction, there is a natural way to identify two disks which are
already attached along their boundaries.  Once we perform these
identifications, we next look at the 3-skeleton.  We now have eight
copies of Figure~\ref{fig-solid} attached to the 2-skeleton, one for
each of the distinct label-preserving automorphisms of its boundary.
Again, we ignore the eight distinct ways in which these are mapped to
the Poincar\'{e} construction and simply identify the interiors.  The
result in this case is the structure shown in Figure~\ref{fig-solid}.  In
general the Cayley category is infinite.  This example is a special
case since the relator contains the entire Cayley graph in its
1-skeleton.

\end{example}

The height of a general relator defines a filtration on any set of
general relators that is sufficient for the above constructions, but
it is a filtration which is often not fine enough for more involved
operations. For this reason a rank function is allowed which provides
an alternative way to order the general relators in a presentation.

\begin{define}\label{def-rank}

{\bf [Rank of a Relator]} A rank function on a presentation $G = \langle
A | \relators \rangle$ is a function from $A \cup \relators$ to
the natural numbers which assigns each generator a rank of 1, and all
of the general relators in $\relators$ a rank of at least 2.  The only
restriction is that if $R$ is a general relator which is contained in
the boundary of the general relator $S$, then $\rank(R) < \rank(S)$. 

\end{define}

\section{General Small Cancellation Theory}\label{sect-gsct}

A very brief overview of general small cancellation theory will now be
given.  The reader is referred to \cite{mccammond-gsct} for further
details and for the definition of any concepts which are not defined
here. A general small cancellation presentation of a group $G$ is a
measured presentation $G = \langle A | \relators \rangle$ which
satisfies the axioms of general small cancellation theory.  The
presentation is measured in the sense that each of the general relators
in the set $\relators$ is required to possess a function which measures
the length of paths in the boundary of the relator for the purposes of
the small cancellation conditions.  Following a discussion of these
metrics, the axioms of general small cancellation theory will be
discussed in detail.  The section concludes with a statement of the
results derived in \cite{mccammond-gsct}.

\subsection{Relator Metrics}

To generalize the small cancellation hypotheses to the context of
general relators requires the introduction of functions which measure
the extent to which a particular path wraps around the boundary of a
general relator.  The function on paths in $R$ will be called the
relator metric for $R$ and it will be denoted $d_R$.  The use of the
term `metric' is justified by the fact that the function $d_R$ will
induce a metric on the points in the universal cover of the boundary
of $R$.  The exact properties required of the function $d_R$ are
listed below.

\begin{define}\label{def-relatormetrics}

{\bf [Relator Metrics]}  A relator metric for a general relator $R$ is a
function $d_R$ which assigns a nonnegative real number to every path
in $R$.  The function $d_R$ must also satisfy the following six
properties:

\begin{enumerate}

\item $d_R(U) = d_R(V)$ whenever $UV^{-1}$ is a contractible loop in
$\partial R$

\item $d_R(U) \geq 0$ and $d_R(U) = 0$ iff $U$ is a contractible loop
in $\partial R$

\item $d_R(U) = d_R(U^{-1})$

\item $d_R(UV) \leq d_R(U) + d_R(V)$

\item if $U$ is a path which forms a loop in $\partial R$ then $d_R(U)
\geq \windnum(U)$

\item if $U$ and $V$ are paths in $R$ which differ by an automorphism
of $R$ then $d_R(U) = d_R(V)$.

\end{enumerate}

\end{define}

One example of a relator metric which is always available is the
normalized graph metric defined in the previous section. Although the
graph metric always satisfies the definition, the added flexibility of
allowing other metrics will be retained since other possibilities prove
to be useful in the construction of the general small cancellation
presentation for the Burnside groups.  A general relator together with
a specified relator metric is said to be measured.  Similarly, a
general presentation $G = \langle A | \relators \rangle$ over a set of
measured relators is called a measured presentation.

\begin{define}\label{def-reducedwords}

{\bf [Reduced Words]} Corresponding to the various ways of measuring
the length of a path, there are several ways of describing the
`tautness' of a word $W$.  Most familiarly, a geodesic is a word which
is not equivalent to any strictly shorter word in $G$.  The word $W$
will be called Dehn-reduced if there do not exist words $U$ and $V$
such that $U$ is a subword of $W$, $UV$ is readable as a loop in one
of the general relators in $\relators$, and the length of $V$ is
strictly less than the length of $V$.  As in the traditional small
cancellation theory, a geodesic is always Dehn-reduced but the
converse is false.  Finally, $W$ is said to be $\mu$-free if it is
reduced in the free group and does not contain more than $\mu$ of a
relator in $\relators$ as measured by its relator metric.
Specifically, there cannot exist a word $U$ and a relator $R$ such
that $U$ is a subword of $W$ as well as a path in $R$ and $d_R(U) >
\mu$.

\end{define}

The multiplicity of metrics and the variety of corresponding reductions
arise from various facets implicit in small cancellation theory.  For
instance, the curvature condition on van Kampen diagrams arises from
the purely combinatorial and decidedly nonmetric condition $C(6)$,
while the strict shortening of the length of a word is necessary in
many cases to guarantee that a process will stop in a finite number of
steps.  As long as the metrics are suitably related to each other,
different metrics can fulfill these two roles. In the present system
these roles will be performed by the relator metrics and the
(normalized) graph metrics, respectively.

\subsection{The Axioms}

This particular version of general cancellation theory requires seven
axioms involving five constants. While the numbers may seem large, the
axioms are far from independent, as are the constants.  The
redundancy built into the system is intended to improve the overall
conceptual clarity.  After the constants and the axioms
are described informally, the exact statements will be given.

The five constants ($\alpha, \beta, \gamma, \delta, \textrm{ and }
\epsilon$) are functionally defined in the first five axioms.  Of
these, the constant $\alpha$ is the most analogous to the constant
used in traditional small cancellation theory in that it measures the
degree to which one relator can be contained in another without being
subsumed.  Specifically, the first axiom states that if $U$ is a path
which can be found in two distinct general relators (or in one general
relator in two non-isomorphic ways) then either the measure of $U$ in
the relator metrics is small (less than $\alpha$) or else the general
relator $R$ for which $d_R(U) \geq \alpha$ is contained inside the
boundary of the other relator.  This can be said more precisely using
maps.  The second axiom states that if a lower ranked relator overlaps
with a higher ranked relator then any path contained in the overlap is
small ($< \beta$) when measured by the relator metric for the higher
ranked relator.  The situation described in the second axiom is in
some sense redundant since whenever the first axiom holds for some
$\alpha$, the second axiom holds for $\beta$ equal to the same
constant.

Axiom three specifies that the ratio of the width of a general relator
to its length must be at most $\gamma$.  The fourth states that the
value given by the relator metric $d_R$ never differs from the value
given by the normalized graph metric $| \cdot |_R$ by more than
$\delta$. The fifth is probably the least used.  If $U$ is the
shortest path from a vertex in $R$ to a loop in $R$ with a non-trivial
winding number, then the measure of $U$ using the relator metric is at
most $\epsilon$. The seventh axiom describes a few convenient
relationships between the constants described above.  The only
constraints which have been included are those needed to derive the
most basic results for general small cancellation theory.  Other
restrictions on the constants will be listed in the statements of the
theorems as they are needed.

The only axiom which has not yet been described is Axiom~\ref{self}.
There is a situation which arises in the reduction of a van Kampen
diagram to a reduced form in which a 2-cell in the diagram overlaps
with itself to form an annular subdiagram.  In traditional small
cancellation theory the overlapping path can contain at most a subword
of a piece since anything longer would imply that a word in the free
group is conjugate in the free group to its own inverse, and this is a
contradiction.  In the case of a representative loop of a general
relator, however, such a long overlap is possible, and the sixth axiom
simply states that should this situation arise, there is always a way
to reduce the diagram further. Figure~\ref{fig-self} shows an example
of how this can occur.

\begin{define}\label{def-axioms}

{\bf [The Axioms]}  The axioms of general small cancellation theory are
as follows:

\setcounter{axiom}{0}

\begin{axiom}\label{aclosed}

There is a constant $\alpha$ such that every general relator $R \in
\relators$ is $\alpha$-closed with respect to $\relators$.  In
particular, if $U$ is a word readable in general relators $R$ and $S$
via $\relators$-functors $f$ and $g$ respectively, and $d_R(U) \geq
\alpha$, then, since $S$ is $\alpha$-closed, there exists a unique
$\relators$-functor $h:R \to S$ such that $hf = g$.

\end{axiom}

\begin{axiom}\label{upbound}

There is a constant $\beta$ such that whenever a word $U$ is readable
in general relators $R, S \in \relators$ by $\relators$-functors $f$
and $g$ respectively, and either $\rank(R) < \rank(S)$ or $\rank(R) =
\rank(S)$ but there does not exist an $\relators$-functor $h: R \to S$
with $hf = g$, then $d_S(U) < \beta$.

\end{axiom}

\begin{axiom}\label{thin}

There is a constant $\gamma$ such that for all general relators $R \in
\relators$, $\omega_R \leq \gamma |R|$.

\end{axiom}

\begin{axiom}\label{quasi}

There is a constant $\delta$ such that the length of a path $U$ in a
general relator $R \in \relators$ in the relator metric $d_R$ is within
$\delta$ of its length in the normalized graph metric on the boundary
of $R$.  Specifically $| d_R(U) - |U|_R | \leq \delta$.

\end{axiom}

\begin{axiom}\label{zimbound}

There is a constant $\epsilon$ such that whenever $U$ is the shortest
possible path from a vertex in a general relator $R \in \relators$ to a
loop with nonzero winding number in $R$, the length of $U$ in the
relator metric is at most $\epsilon$.  That is, $d_R(U) \leq \epsilon$.

\end{axiom}

\begin{axiom}\label{self}

If $W = XUYU^{-1}$ is a representative of a general relator $R$ in which
both instances of $U$ are properly oriented with respect to $W$, and
$d_R(U) \geq \alpha$, then there exists a word  $V$ such that the cycle
$XVYV^{-1}$ is readable as a contractible loop in $\partial R$
extending the reading of $X$ given by $W$.  The cycle thus bounds a
connected and simply connected $\relators$-diagram $\Delta$ with
$\rank(\Delta) < \rank(R)$.

\end{axiom}

\begin{axiom}\label{constraints}

The constants $\alpha$, $\beta$, $\gamma$, $\delta$, and $\epsilon$
satisfy the following constraints: $\beta \leq \alpha$, and $\gamma,
\delta, \epsilon < \alpha$, and $2 \gamma + \delta \leq \alpha \leq
\frac{1}{6}$.

\end{axiom}

\end{define}

A general small cancellation presentation $G = \langle  A | \relators
\rangle$ is a measured presentation which satisfies the axioms listed
above.  A group which possesses a general small cancellation
presentation is called a general small cancellation group. The key
example of an interesting general small cancellation group is, of
course, the Burnside groups. In \cite{mccammond-burnside}, a general
small cancellation presentation is constructed for each of the free
Burnside groups with exponent greater than or equal to 1260. The
details of the construction of this general presentation are discussed in 
section~\ref{sect-burnside}.

\subsection{The Results}

Once the proper definitions are in place the proofs of the usual small
cancellation results are slight variations on the traditional proofs.
It thus follows quickly that subject to a few restrictions on the
values of the constants, the word problem and the conjugacy problem
are decidable.

\begin{thm}[\cite{mccammond-gsct}]\label{thm-wordprob}	% Lemma 10.21

If $G = \langle A | \relators \rangle$ is a general small cancellation
presentation with $3 \alpha + 2 \gamma + \delta \leq \frac{1}{2}$, then
the group $G$ has a decidable word problem. In particular, these
results are true when $\alpha \leq \frac{1}{6}$ and $\gamma = \delta =
0$ or whenever $\alpha \leq \frac{1}{8}$.

\end{thm}

\begin{thm}[\cite{mccammond-gsct}]\label{thm-conjprob}  % Lemma 10.23

If $G = \langle A | \relators \rangle$ is a general small cancellation
presentation with $4 \alpha + 2 \gamma + \delta \leq \frac{1}{2}$, then
the group $G$ has a decidable conjugacy problem. In particular, this is
true when $\alpha \leq \frac{1}{8}$ and $\gamma = \delta = 0$ or
whenever $\alpha \leq \frac{1}{10}$.

\end{thm}

Similarly, a finitely presented general small cancellation presentation
is word-hyperbolic.  The proof of this particular theorem uses small
cancellation conditions to prove that these groups satisfy Rip's thin
triangle conditon.

\begin{thm}[\cite{mccammond-gsct}]\label{thm-hyperbolic}  % Lemma 11.19

If $G = \langle A | \relators \rangle$ is a finitely presented general
small cancellation presentation with $3 \alpha + \delta \leq
\frac{1}{3}$ and $\alpha \leq \frac{1}{10}$, then the group $G$ is
word hyperbolic.  In particular, the result is true when $\alpha =
\frac{1}{10}$ and $\gamma = \delta = 0$, or whenever $\alpha \leq
\frac{1}{12}$.

\end{thm}

In the traditional theory, the Cayley complex of a small
cancellation presentation is shown to be aspherical.  The following is
the general small cancellation version of this fact.

\begin{thm}[\cite{mccammond-gsct}]\label{thm-contractible}  % Lemma 12.12

If $G = \langle A | \relators \rangle$ is a general small cancellation
presentation, and $C$ is the Cayley category of the presentation, then
$C$ is contractible.

\end{thm}

As discussed earlier, the difference between the universal cover of the
Poincar\'{e} construction and the Cayley category of a general
presentation is determined by the automorphism groups of the relators.
As a result, presentations whose relators have no nontrivial
automorphisms have additional properties.

\begin{thm}[\cite{mccammond-gsct}]\label{thm-torsion}  % Lemma 12.14

If $G = \langle A | \relators \rangle$ is a finitely presented general
small cancellation presentation with $\alpha \leq \frac{1}{10}$, then
the following five conditions are equivalent:

(1) the group $G$ is torsion-free

(2) all of the general relators in $\relators$ have no nontrivial
automorphisms

(3) the universal cover of the Poincar\'{e} construction is collapsed

(4) the universal cover of the Poincar\'{e} construction is
contractible

(5) the Poincar\'{e} construction is a $\textrm{K}(G,1)$-space.

\end{thm}

The final result of general small cancellation theory concerns the
relationship between the automorphism groups of the general relators in
a presentation and the finite subgroups of the group described by the
presentation.  Under fairly mild restrictions the two lists are identical.

\begin{thm}[\cite{mccammond-gsct}]	% Lemma 14.17

If $G = \langle A | \relators \rangle$ is a general small cancellation
presentation with $2 \beta + 2 \gamma + \delta \leq \alpha \leq
\frac{1}{12}$ in which $\str(W)$ is finite and effectively
constructible for all words $W \in A^*$, and such that all general
relators in $\relators$ have at least one crucial cone in their
boundary, then every finite subgroup of $G$ is a subgroup of the
automorphism group of some general relator in $\relators$.

\end{thm}

Since the automorphism groups of general relators possessing crucial
cones were described in Theorem~\ref{thm-autR}, the following corollary
is immediate.

\begin{cor}[\cite{mccammond-gsct}]

If $G = \langle A | \relators \rangle$ is a general small cancellation
presentation with $2 \beta + 2 \gamma + \delta \leq \alpha \leq
\frac{1}{12}$ in which $\str(W)$ is finite and effectively
constructible for all words $W \in A^*$, and such that all general
relators in $\relators$ have at least one crucial cone in their
boundary, then each finite subgroup of $G$ can be embedded in the
extension of a dihedral group by a 2-group. As a consequence, a general
small cancellation group (satisfying the above conditions) which does
not contain any elements of even order has only cyclic finite
subgroups.

\end{cor}

\section{Burnside Groups}\label{sect-burnside}

Before discussing the results concerning the Burnside groups in
greater detail, the reader is reminded that the descriptions given in
this final section are of a preliminary nature.  In particular, the
exponent 1260 is especially provisional.

In the tradition of Adian, Novikov, Ol'shanskii, and Ivanov, the
general small cancellation presentation of the Burnside groups,
$\burnside(m,n) = \langle A | \relators \rangle$, is constructed
inductively. In particular, most of the lemmas in
\cite{mccammond-burnside} are proved by simultaneous induction on a
single parameter $k$.  The parameter $k$ refers to the ranks of the
general relators used in the presentation.  The notations used are as
follows: $\relators_k$ represents the set of general relators of rank
$k$, $\relators(k)$ represents the set of general relators of rank at
most $k$, and $\relators$ is the union of the sets $\relators_k$ over
all positive integers $k$.  There is a corresponding list of groups:
$G(k) = \langle A | \relators(k) \rangle$, and $G = \langle A |
\relators \rangle$.  After a few of the details of the inductive
construction are presented, the conclusions which result from the
existence of a general small cancellation presentation for the
Burnside groups will be listed.  We begin with a simple illustration
of the approach.

\begin{example}

{\bf [The Main Idea]} Consider a presentation of the form $G = \langle
A | X^n, Y^n \rangle$ where $n$ is large and $X$ and $Y$ are simple
words.  It is well-known that the only way that a high power of $X$
can occur as a subword of $Y^n$ is if the length of $X$ is small
compared to that of $Y$.  More precisely, if $X^i$ is a subword of
$Y^n$ then it is also a subword of a conjugate of $Y^2$.  In this
case, the construction described in outline below will attach disks
with boundaries labeled $X^n$ to each of the subwords $X^i$ in
$Y^n$.  The result will be somewhat similar to the general relator
shown in Figure~\ref{fig-solid} except that that general relator
corresponds to the word $(abcb)^2$ instead of a higher power.  In
general, we will begin with a word $Y^n$ and then we will attach
previously constructed general relators whenever they have a high
power of their defining word contained in the general relator under
construction.  Under appropriate restrictions, this process of
attaching general relators stops, and the result has sufficiently nice
properties to be of use in an inductive construction.

\end{example}

\subsection{The Inductive Construction}

\begin{figure}
	\begin{center}~
		\epsffile{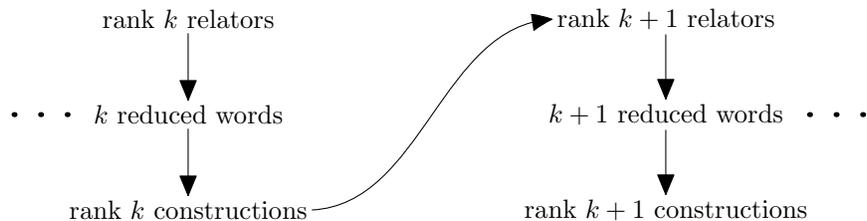}
	\end{center}
	\caption{The order of definition \label{fig-order}}
\end{figure}

The general relators used to present the Burnside groups are defined one
rank at a time.  The rank, $k$, indexes not only a set of general
relators, but also a set of words, a set of cycles, and a set of
constructions.  The definitions of these concepts are intertwined and
inductively defined.  The cycle of definitions is as follows: reduced
words and cycles in rank $k$ are used to construct the rank $k$
straightline and circular constructions.  These constructions are then
used to define the rank $k+1$ general relators, which are in turn used
to define reduced words and cycles in rank $k+1$. See
Figure~\ref{fig-order}.

Before describing the induction itself in more detail, a word should be
said about the values of the constants used in the inductive construction
and about the definitions in low cases.  The values of the constants used 
are as follows:

\[ 
	\alpha = \frac{1}{12}, 
	\beta = \frac{1}{210}, 
	\gamma = \frac{1}{70},
	\delta = \frac{1}{30}, 
	\epsilon = \frac{1}{630}, 
	n \geq 1260 
\]

\noindent These values are sufficient to satisfy all of the additional
restrictions placed on the value of the constants in the statements of
the results listed in Section~\ref{sect-gsct}.

The induction starts as follows.  Since by definition all general
relators have a rank of at least 2, there are no rank 1 relators and
$\relators_1 = \emptyset$.  The cycle of definitions thus really begins
with the notions of 1-reduced words and 1-reduced cycles. A word
(cycle) is called 1-reduced iff it is reduced (cyclically reduced) in
the free group.  If $W$ is a 1-reduced word, then $\str_1(W)$ is simply
the abstract path labeled by the word $W$. Similarly, if $W$ is a
1-reduced cycle then $\cir_1(W)$ is the abstract loop labeled by the
cycle $W$. The definitions are then extended to arbitrary words
(respectively, cycles) which are not 1-reduced by defining the rank 1
straightline (circular) construction on the word (cycle) to be the
appropriate construction of its reduction in the free group.  The rank
2 general relators are defined from the rank 1 constructions. In
particular, if $W$ is a simple word in the traditional sense (i.e. not
equal to a proper power of a shorter word), and if the cycle of $W$
does not contain $\beta n$ powers of any simple word, then the cone
over the construction $\cir_1(W^n)$ is a rank 2 general relator.  The
set of all rank 2 general relators is called $\relators_2$.  Since
these relators do not contain any $\beta n$ powers of any simple word,
they will  not contain $\beta$ of the boundary of any of the other rank
2 relators. Thus the rank 2 relators satisfy the traditional small
cancellation hypothesis $C'(\beta)$. As a consequence, the group $G(2)$
is a traditional small cancellation group.

The rank 1 straightline and circular constructions have already been
defined.  The rank 2 versions of these constructions will be
illustrated using Figure~\ref{fig-moebius}.  Let $G = \langle a,b,c |
(ab)^3 = (bc)^3 = (ca)^3 \rangle$, let $W =babcba^{-1}c^{-1}$, and let
the three relators in $G$ be assigned a rank of 2.  The construction
$\str_2(W)$ is shown in Figure~\ref{fig-moebius}.  The construction
$\cir_2(W)$ is the structure obtained by gluing the left and the right
edges according to orientation.

\begin{define}\label{def-str}

{\bf [The Straightline Construction]}  The rank $k$ straightline
construction on a $k$-reduced word $W$ is in some sense the smallest
structure which contains $W$ as a path and which is $\alpha$-closed with
respect to $\relators(k)$.  The latter condition means that given any
path $U$ in the structure and a general relator $R \in \relators(k)$ such
that $d_R(U) \geq \alpha$, the relator $R$ is already attached to the
structure along this path $U$.  

\end{define}

\begin{define}\label{def-cir}

{\bf [The Circular Construction]}  The circular construction is defined
similarly.  The rank $k$ circular construction on a $k$-reduced cycle
$W$ is the smallest structure which contains $W$ as a loop and is
$\alpha$-closed with respect to $\relators(k)$.

\end{define}

If $G = \langle A | \relators \rangle$ is any general small cancellation
presentation with $\alpha \leq \frac{1}{8}$, it is shown in
\cite{mccammond-gsct} that the rank $k$ straightline and circular
constructions on $W$ exist and are well-defined. Under certain
additional restrictions, the constructions which result are finite. 
The presentation constructed inductively for the Burnside groups
satisfies these restrictions, and it is the rank $k$ circular
constructions which are used to define the boundaries of what will become
the rank $k+1$ general relators.  Once these relators have been
constructed, the axioms of general small cancellation are shown to hold. 
At that point the results listed in the previous section become
available, and the induction continues.

\subsection{The Results}

The most important result on the Burnside groups of sufficiently large
exponent that is established in \cite{mccammond-burnside} is that they
possess a general small cancellation presentation.

\begin{thm}\label{thm-burnside=gsc}

For $n \geq 1260$ the Burnside group $\burnside(m,n)$ possesses 
a general small cancellation presentation.

\end{thm}

Most of the other results for these groups follow immediately from the
existence of such a presentation, and from the general small
cancellation theory described in the previous section.  In particular,
the following corollary is immediate from Theorem~\ref{thm-wordprob}
and Theorem~\ref{thm-conjprob}.

\begin{cor}

For all $n \geq 1260$, the Burnside group $\burnside(m,n)$ has a decidable
word problem and a decidable conjugacy problem.

\end{cor}

As might be expected, the classification of the structure of the finite
subgroups of the Burnside groups is more specific than the structure of
the finite subgroups in a general small cancellation group in general.

\begin{thm}

Every finite subgroup of the group $\burnside(m,n)$ is contained in a
direct product of a dihedral of order $2n$ with a finite number of
dihedral $2$-groups whose exponent divides $n$.  

\end{thm}

It is from the classification of the finite subgroups that the last two
major results are derived.

\begin{cor}

For all $n \geq 1260$ and all $m \geq 2$, the Burnside group
$\burnside(m,n)$ is infinite.

\end{cor}

\begin{proof}

If $\burnside(m,n)$ were finite, then it would itself be a finite
subgroup.   Thus, once it is shown that $\burnside(m,n)$ cannot be embedded in a
finite product of dihedral groups, the proof will be complete.  It is an
easy matter to demonstrate this.

\end{proof}

\begin{cor}

For all $n \geq 1260$ and all $m \geq 2$, the Burnside group
$\burnside(m,n)$ is not finitely presented.

\end{cor}

\begin{proof}

If $\burnside(m,n)$ were finitely presented, then by
Theorem~\ref{thm-hyperbolic} it would be a word hyperbolic group, and
thus automatic.  However, it is well known that finitely presented
infinite torsion automatic groups do not exist.

\end{proof}


\newpage

\begin{thebibliography}{1}

\bibitem{adian-english}
S.~Adian.
\newblock {\em The {B}urnside Problem and Identities in Groups}.
\newblock Springer-Verlag, New York, 1979.

\bibitem{burnside}
W.~Burnside.
\newblock {\em An unsettled question in the theory of discontinuous groups.}
\newblock Quart. J. of Pure and Appl. Math., {\bf 33} (1902), 230--238.

\bibitem{ivanov}
S.~Ivanov.
\newblock {\em The free {B}urnside groups of sufficiently large exponents.}
\newblock International Journal of Algebra and Computation, {\bf 4} (1994),
1--308.

\bibitem{lysionok}
I.~Lysionok.
\newblock {\em Infinite {B}urnside groups of even period.}
\newblock Izv. Ross. Akad. Nauk Ser. Mat. {\bf 60} (1996), 3--224.

\bibitem{mccammond-gsct} 
J.~McCammond.  
\newblock {\em A general small cancellation theory} 
\newblock to appear in Int. J. of Alg. and Comp.

\bibitem{mccammond-burnside}
J.~McCammond.
\newblock {\em The {B}urnside groups via a general small cancellation theory.}
\newblock In preparation.

\bibitem{adiannovikov}
P.~Novikov and S.~Adian.
\newblock {\em On infinite periodic groups, {I}, {II}, {III}.}
\newblock Izv. Acad. Nauk. SSSR Ser. Matem.,
  {\bf 32} (1968), 212--244,251--524,709--731.

\bibitem{olshanskii}
A.~Yu. Ol'shanskii.
\newblock {\em On the {N}ovikov-{A}dian theorem.}
\newblock Math. Sbornik, {\bf 118} (1982), 203--235.

\bibitem{quillen}
D.~Quillen.
\newblock {\em Higher Algebraic K-Theory: I}. 
\newblock Springer-Verlag, New York, 1973.

\bibitem{rips}
E.~Rips.
\newblock {\em Generalized small cancellation theory and applications {I}.}
\newblock Isr. J. Math. {\bf 41} (1982), 1--146.

\end{thebibliography}
\end{document}